# Density in the half-line Schwartz space of functions whose Fourier-Laplace transform has natural boundary the real line

Andreas Chatziafratis [1,*] and Telemachos Hatziafratis [2]

**Abstract.** In this note, we prove that the Fourier-Laplace transform of the "typical" function (i.e., generic in the sense of Baire's category theorem) in the Schwartz class of the half-line, being analytic in the lower half of the complex plane, has natural boundary the axis of the real numbers. We also provide variations and generalizations.

## 1. Introduction

In certain problems of classical applied analysis, mathematical physics and partial differential equations (e.g. addressing initial-boundary-value problems posed on quarter-planes; see, for example, the recent line of investigations [1-17] and many more works [18-61]), it is customary to deal with functions (e.g. initial data) that are defined on the half-line and belong to the corresponding Schwartz class (or an appropriate class of "sufficient" smoothness and decay). The "half" Fourier transform (HFT) of such functions is then almost invariably employed during the analysis of such problems. In this case, it is easy to see that the HFT is well-defined in the lower-half complex plane (LHP) since the assumed decay guarantees absolute convergence for the associated integral (and some of its derivatives). In particular for Schwartz functions, the HFT defines an analytic function in the LHP which is infinitely smooth up to the real axis. On the other hand, it is a trivial task to find specific examples of functions belonging to the Schwartz class whose HFT extends to an analytic function in parts of the upper half plane too, or even to an entire holomorphic function. A reasonable question thus arises: *Is there a Schwartz function, supported on the positive half-line, whose Fourier transform does not extend holomorphically beyond the real axis?* We herein settle this matter by first constructing functions which have singularities at any given finite set of points on the real line, and then appealing to a suitable implementation of the well-known Baire category theorem in order to prove that there exist such functions whose half Fourier transform's *natural* boundary is the real axis. (For a discussion of natural boundaries, see e.g. [62].) As a matter of fact, we show that "such" functions are "dense" in the Schwartz class with respect to its usual topology. This result may not be entirely surprising, nonetheless it appears to be missing from the literature. It may additionally be regarded as a small contribution towards the description of the range of half Fourier-Laplace transform. (For discussions on this transform, Paley-Wiener-type theorems and related topics, we refer to, e.g., [63-66].) Furthermore, we extend the results to broader function classes and provide variations too.

As we pointed out, the key tool, in proving the density – and therefore the existence – of such functions, is the fundamental Baire category theorem. It is noted that this theorem has been invoked in numerous proofs of existence and density theorems: indicatively, see for example [67-72]; some further applications and historical comments may be found, for instance, in the review paper [73].

*Notation and preliminaries* Let $\mathcal{S}([0,\infty))$ be the space of Schwartz functions in the half-line. More precisely,

$$\mathcal{S}([0,\infty)) := \{f \in C^\infty([0,\infty)) : \rho_\ell(f) < \infty \text{ for every nonnegative integer } \ell\}$$

where

$$\rho_\ell(f) := \sup\left\{(1+t)^\ell \left|\frac{\partial^k f(t)}{\partial t^k}\right| : k \leq \ell, t \geq 0\right\}.$$

---

[1] Department of Mathematics, School of Science, National and Kapodistrian University of Athens, Greece,
  and Institute of Applied & Computational Mathematics, FORTH, Crete, Greece

[2] Department of Mathematics, School of Science, National and Kapodistrian University of Athens, Greece

*corresponding author, e-mail: chatziafrati@math.uoa.gr

Endowed with the topology defined by the norms $\{\rho_\ell : \ell\}$, $\mathcal{S}([0,\infty))$ becomes a Fréchet space. We recall that this means that the topology of $\mathcal{S}([0,\infty))$ is defined by the metric

$$\rho(f,g) := \sum_\ell \frac{1}{2^\ell} \frac{\rho_\ell(f-g)}{1+\rho_\ell(f-g)}, \quad f,g \in \mathcal{S}([0,\infty)).$$

Then $\mathcal{S}([0,\infty))$ becomes a complete metric space and a sequence $f_j \to g$, in $\mathcal{S}([0,\infty))$ with respect to the metric $\rho$, as $j \to \infty$, if and only if $\rho_\ell(f_j - g) \to 0$, as $j \to \infty$, for every $\ell$.

For $f \in \mathcal{S}([0,\infty))$, we consider its Fourier-Laplace tranform

$$\hat{f}(z) = \int_0^\infty e^{-izt} f(t) dt, \text{ defined for } z \in \mathbb{C} \text{ with } \operatorname{Im} z \leq 0.$$

Then the function $\hat{f}(z)$ is $C^\infty$ in $\{z \in \mathbb{C} : \operatorname{Im} z \leq 0\}$ and holomorphic (analytic) for $z \in \mathbb{C}$ with $\operatorname{Im} z < 0$. This follows from Lebesgue's dominated convergence theorem and the fact that $\sup\{t^\ell |f(t)| : t \geq 0\} < \infty$, $\forall \ell$.

Also, integration by parts gives

$$\hat{f}(z) = \frac{f(0)}{iz} + \frac{f'(0)}{(iz)^2} + \cdots + \frac{f^{(k-1)}(0)}{(iz)^k} + \frac{1}{(iz)^k} \int_0^\infty e^{-izt} f^{(k)}(t) dt, \text{ for } z \in \mathbb{C} - \{0\} \text{ with } \operatorname{Im} z \leq 0,$$

and, in particular, we have

$$\hat{f}(z) = \frac{f(0)}{iz} + \frac{f'(0)}{(iz)^2} + \cdots + \frac{f^{(k-1)}(0)}{(iz)^k} + O(1/z^{k+1}), \text{ as } z \to \infty, \text{ with } z \in \mathbb{C}, \operatorname{Im} z \leq 0.$$

For example, if $f(t) = e^{-t^p}$, $t \geq 0$, for some $p > 1$, the function $\hat{f}(z)$ extends to an entire analytic function, while if $f(t) = t^\nu e^{-\sigma t}$, $t \geq 0$, for some $\nu \in \mathbb{N} \cup \{0\}$ and $\sigma \in \mathbb{C}$ with $\operatorname{Re} \sigma > 0$, then $\hat{f}(z) = (-i)^\nu \nu!/(z - i\sigma)^{\nu+1}$, i.e., $\hat{f}(z)$ extends to an analytic function in $\mathbb{C} - \{i\sigma\}$ with a pole at $z = i\sigma$, of order $\nu + 1$.

In this note we will show that "in general" the function $\hat{f}(z)$, being holomorphic in $\{z \in \mathbb{C} : \operatorname{Im} z < 0\}$, has natural boundary the real line $\{z \in \mathbb{C} : \operatorname{Im} z = 0\}$, i.e., it does not extend analytically beyond any point $\alpha \in \mathbb{R}$. Let us point out that $\hat{f}(z)$ does not extend analytically beyond some point $\alpha \in \mathbb{R}$ if and only if

$$r(f,\alpha) := \limsup_{n \to \infty} \sqrt[n]{\frac{|(\hat{f})^{(n)}(\alpha)|}{n!}} = \infty.$$

As we mentioned before, the restriction $\hat{f}|_{\mathbb{R}}$, of $\hat{f}(z)$ to the real line, is $C^\infty$, and the numbers $(\hat{f})^{(n)}(\alpha)$ are the real derivatives of the function $\hat{f}(x)$, $x \in \mathbb{R}$, evaluated at $x = \alpha$.

**Notation** For $f \in \mathcal{S}([0,\infty))$, $\alpha \in \mathbb{R}$ and $n \in \mathbb{N}$, we write

$$C(f,\alpha,n) := \sqrt[n]{\frac{|(\hat{f})^{(n)}(\alpha)|}{n!}}.$$

More precisely, we will prove the following theorem.

**Theorem 1** *The set $\{f \in \mathcal{S}([0,\infty)) : r(f,\alpha) = \infty \text{ for every } \alpha \in \mathbb{R}\}$ is dense and $\mathcal{G}_\delta$ in the space $\mathcal{S}([0,\infty))$.*

## 2. Proof of Theorem 1

We start with an example of a function in $\mathcal{S}([0,\infty))$, whose Fourier transform has a singular point on the real line. We will need such an example in the proof of Theorem 1.

**Example** For $t \geq 0$, we define the function $\chi_0(t) := e^{-\sqrt{t}}$. Its Fourier transform is

$$\hat{\chi}_0(z) = \int_0^\infty e^{-izt} e^{-\sqrt{t}} dt = 2\int_0^\infty e^{-izs^2} e^{-s} s\, ds, \quad z \in \mathbb{C}, \; \text{Im}\, z \leq 0.$$

Computing the derivatives

$$\hat{\chi}_0^{(n)}(z) = 2\int_0^\infty e^{-izs^2} e^{-s} (-is^2)^n s\, ds,$$

we find that

$$\hat{\chi}_0^{(n)}(0) = 2\int_0^\infty e^{-s}(-is^2)^n s\, ds = 2(-i)^n \int_0^\infty e^{-s} s^{2n+1} ds = 2(-i)^n (2n+1)!.$$

It follows that the radius of convergence of the power series

$$\sum_{n=0}^\infty \frac{\hat{\chi}_0^{(n)}(0)}{n!} z^n$$

is zero. Thus, the function $\hat{\chi}_0(z)$, which is holomorphic in $\{\text{Im}\, z < 0\}$, does not extends analytically beyond the point $0$.

More generally, if $\chi_\alpha(t) = e^{i\alpha t} e^{-\sqrt{t}}$, $\alpha \in \mathbb{R}$, the radius of convergence of the power series

$$\sum_{n=0}^\infty \frac{\hat{\chi}_\alpha^{(n)}(\alpha)}{n!} (z-\alpha)^n$$

is zero and the function $\hat{\chi}_\alpha(z)$, which is holomorphic in $\{\text{Im}\, z < 0\}$, does not extend analytically beyond the point $\alpha$.

**Analytic continuation of $\hat{\chi}_0(z)$ in $\mathbb{C} - i(\mathbb{R}^+ \cup \{0\})$** For this purpose, let us compute the values of the function $\hat{\chi}_0(z)$ at the points $z = -i\mu$ when $\mu > 0$. We have

$$\hat{\chi}_0(-i\mu) = \int_0^\infty e^{-\mu t} e^{-\sqrt{t}} dt = 2\int_0^\infty e^{-\mu s^2} e^{-s} s\, ds$$

$$= \frac{1}{\mu} + \frac{1}{\mu}\int_0^\infty e^{-\mu s^2} e^{-s} ds = \frac{1}{\mu} + \frac{\sqrt{\pi}}{2} \frac{1}{\mu\sqrt{\mu}} \exp\left(\frac{1}{4\mu}\right)\left[1 - \int_0^{1/(2\sqrt{\mu})} e^{-t^2} dt\right]$$

$$= \frac{1}{\mu} + \frac{\sqrt{\pi}}{2} \frac{1}{\mu\sqrt{\mu}} \exp\left(\frac{1}{4\mu}\right) - \frac{\sqrt{\pi}}{2} \frac{1}{\mu\sqrt{\mu}} \exp\left(\frac{1}{4\mu}\right) \int_0^{1/(2\sqrt{\mu})} e^{-t^2} dt$$

$$= \frac{1}{\mu} + \frac{\sqrt{\pi}}{2} \frac{1}{\mu\sqrt{\mu}} \exp\left(\frac{1}{4\mu}\right) - \frac{\sqrt{\pi}}{2} \frac{1}{\mu\sqrt{\mu}} \exp\left(\frac{1}{4\mu}\right) \sum_{n=0}^\infty \frac{(-1)^n}{n!(2n+1)} \left(\frac{1}{2\sqrt{\mu}}\right)^{2n+1}$$

$$= \frac{1}{\mu} + \frac{\sqrt{\pi}}{2} \frac{1}{\mu\sqrt{\mu}} \exp\left(\frac{1}{4\mu}\right) - \frac{\sqrt{\pi}}{4} \frac{1}{\mu^2} \exp\left(\frac{1}{4\mu}\right) \sum_{n=0}^\infty \frac{(-1)^n}{4^n n!(2n+1)} \frac{1}{\mu^n}.$$

Setting $\mu = iz$ in the last expression of $\hat{\chi}_0(-i\mu)$, we find

$$\hat{\chi}_0(z) = \frac{1}{iz} + \frac{\sqrt{\pi}}{2} \frac{1}{iz\sqrt{iz}} \exp\left(\frac{1}{4iz}\right) + \frac{\sqrt{\pi}}{4} \frac{1}{z^2} \exp\left(\frac{1}{4iz}\right) \sum_{n=0}^\infty \frac{(-1)^n}{4^n n!(2n+1)} \frac{1}{(iz)^n}.$$

The RHS of the above equality defines the analytic continuation of $\hat{\chi}_0(z)$, for $z \in \mathbb{C} - i(\mathbb{R}^+ \cup \{0\})$, choosing the following branch of the square root: If $z = |z|e^{i\arg z}$, with $-\pi < \arg z \leq \pi$, then $\sqrt{z} = \sqrt{|z|} e^{\arg z/2}$.

**Modification of** $\chi_\alpha$. Since $\lim_{t \to 0^+}[d\chi_\alpha(t)/dt] = \infty$, $\chi_\alpha \notin C^1([0,\infty))$ and, therefore, $\chi_\alpha \notin \mathcal{S}([0,\infty))$. However, we can easily modify the function $\chi_\alpha$ in order to have a function in $\mathcal{S}([0,\infty))$ with similar properties. Indeed, it suffices – for example – to set $\varphi_\alpha(t) = e^{i\alpha t}e^{-\sqrt{t+1}}$, $\alpha \in \mathbb{R}$. Then $\varphi_\alpha \in \mathcal{S}([0,\infty))$ and

$$\hat{\varphi}_\alpha(z) = \int_0^\infty e^{-izt}e^{i\alpha t}e^{-\sqrt{t+1}}dt = e^{iz}e^{-i\alpha}\int_1^\infty e^{-izt}e^{i\alpha t}e^{-\sqrt{t}}dt = e^{iz}e^{-i\alpha}\hat{\chi}_\alpha(z) - e^{iz}e^{-i\alpha}\int_0^1 e^{-izt}e^{i\alpha t}e^{-\sqrt{t}}dt.$$

Since the function $\int_0^1 e^{-izt}e^{i\alpha t}e^{-\sqrt{t}}dt$ is an entire function of $z$, $\hat{\varphi}_\alpha(z)$, which is – a priori – holomorphic in $\{\operatorname{Im} z < 0\}$, extends analytically to $\mathbb{C} - [\alpha + i(\mathbb{R}^+ \cup \{0\})]$, but does not extend analytically beyond the point $\alpha$ and the radius of convergence of the power series $\sum_{n=0}^\infty [\hat{\varphi}_\alpha^{(n)}(\alpha)/n!](z-\alpha)^n$ is zero.

**Remarks** 1 If $f \in \mathcal{S}([0,\infty))$ and $r(f,\alpha) = \infty$ for every $\alpha \in \mathbb{R}$, then the function $\hat{f}\big|_\mathbb{R} \in C^\infty(\mathbb{R})$ and nowhere real analytic, i.e., at no point $x \in \mathbb{R}$.

2 The radius of convergence of the power series $\sum_{n=0}^\infty [\hat{\varphi}_\alpha^{(n)}(\alpha)/n!](z-\alpha)^n$ is equal to $|\alpha|$, for every $\alpha \in \mathbb{R}$. More generally, if $\Phi_\alpha(z)$ is the analytic continuation of $\hat{\varphi}_\alpha$ in $\mathbb{C} - [\alpha + i(\mathbb{R}^+ \cup \{0\})]$, then

$$\limsup_{n \to \infty} \sqrt[n]{|\Phi_\alpha^{(n)}(\sigma)|/n!} = \frac{1}{|\sigma - \alpha|}, \text{ for } \sigma \in \mathbb{C} - [\alpha + i(\mathbb{R}^+ \cup \{0\})].$$

The following lemmas are also needed in the proof of the theorem.

**Lemma 1** *Let $n \in \mathbb{N} \cup \{0\}$ and $z \in \mathbb{C}$, with $\operatorname{Im} z \leq 0$, be fixed. Then the map*

$$\mathcal{S}([0,\infty)) \to \mathbb{C}, \ f \to (\hat{f})^{(n)}(z), \tag{1}$$

*is continuous (in $f \in \mathcal{S}([0,\infty))$).*

***Proof*** Let $f_j \in \mathcal{S}([0,\infty))$ be a sequence which converges to $0$, as $j \to \infty$. We claim that $(\hat{f}_j)^{(n)}(z) \to 0$. Indeed, since

$$\sup\{t^\ell|f_j(t)| : t \geq 0\} \to 0, \text{ as } j \to \infty, \text{ for } \ell = 0,1,\ldots,n+2,$$

it follows that for every $\varepsilon > 0$ there exists $j_0 = j_0(\varepsilon)$ so that

$$j \geq j_0 \Rightarrow \sup\{(1+t)^{n+2}|f_j(t)| : t \geq 0\} < \varepsilon.$$

Since

$$(\hat{f}_j)^{(n)}(z) = \int_0^\infty e^{-izt}(-it)^n f_j(t)dt,$$

we obtain

$$|(\hat{f}_j)^{(n)}(z)| \leq \int_0^\infty e^{t\operatorname{Im}z}t^n|f_j(t)|dt < \binom{n+2}{n}\varepsilon\int_0^\infty \frac{dt}{(1+t)^2} = \binom{n+2}{n}\varepsilon, \text{ for } j \geq j_0.$$

This completes the proof.

**Lemma 2** *Let $\alpha \in \mathbb{R}$ and $N \in \mathbb{N}$. Then the set*

$$\Omega(\alpha, N) := \bigcup_{n \in \mathbb{N}}\{g \in \mathcal{S}([0,\infty)) : C(g,\alpha,n) > N\}$$

*is open and dense in the space $\mathcal{S}([0,\infty))$.*

*Proof* First, since the map (1) is continuous, the set

$$\{g \in \mathcal{S}([0,\infty)) : C(g,\alpha,n) > N\}$$

is open in $\mathcal{S}([0,\infty))$, for every $n$, and, therefore, $\Omega(\alpha,N)$ is, indeed, open.

In order to prove the density of $\Omega(\alpha,N)$ in $\mathcal{S}([0,\infty))$, let us consider an $f \in \mathcal{S}([0,\infty))$. We claim that there exist elements of the set $\Omega(\alpha,N)$, arbitrarily close to $f$. This is trivial in the case $f \in \Omega(\alpha,N)$. Thus, it suffices to consider $f \in \mathcal{S}([0,\infty)) - \Omega(\alpha,N)$. Then

$$C(f,\alpha,n) \leq N, \text{ for every } n \in \mathbb{N}.$$

It follows that $r(f,\alpha) < \infty$ and, therefore, the power series

$$\sum_{n=0}^{\infty} \frac{(\hat{f})^{(n)}(\alpha)}{n!}(z-\alpha)^n$$

defines an analytic extension of $\hat{f}(z)$ in an open neighborhood of the point $\alpha$.

On the other hand,

$$f + \frac{1}{j}\varphi_\alpha \to f, \text{ in the space } \mathcal{S}([0,\infty)), \text{ as } j \to \infty.$$

But

$$f + \frac{1}{j}\varphi_\alpha \in \Omega(\alpha,N).$$

Indeed, this follows from the fact that $\hat{f}(z)$ extends analytically beyond $\alpha$, whereas $\hat{\varphi}_\alpha$ does not, which implies that $r(f + \frac{1}{j}\varphi_\alpha, \alpha) = \infty$.

This proves the density of $\Omega(\alpha,N)$ in $\mathcal{S}([0,\infty))$.

The following lemma is an immediate consequence of Lemma 2 and Baire's theorem.

**Lemma 3** *Let $\alpha \in \mathbb{R}$. Then the set $\bigcap_{N \in \mathbb{N}} \Omega(\alpha,N)$ is dense and $\mathcal{G}_\delta$ in the space $\mathcal{S}([0,\infty))$.*

**Proof of the Theorem 1** Let $A \subset \mathbb{R}$ be a countable dense subset in the real line. We claim that

$$\{f \in \mathcal{S}([0,\infty)) : r(f,\beta) = \infty \text{ for every } \beta \in \mathbb{R}\} = \bigcap_{\alpha \in A} \bigcap_{N \in \mathbb{N}} \Omega(\alpha,N). \quad (2)$$

Since the inclusion "$\subseteq$", in (2), is trivial, it suffices to prove the other inclusion. For this, let us consider an element $f$ of the set $\bigcap_{\alpha \in A} \bigcap_{N \in \mathbb{N}} \Omega(\alpha,N)$. We claim that $r(f,\beta) = \infty$ for every $\beta \in \mathbb{R}$. Let us fix a $\beta \in \mathbb{R}$ and suppose – to reach contradiction – that the function $\hat{f}(z)$ extends to a holomorphic function in the disc $\Delta := \{z \in \mathbb{C} : |z - \beta| < \varepsilon\}$, for some $\varepsilon > 0$. By the density of $A$ in $\mathbb{R}$, there exists an $\alpha_0 \in A$ so that $|\alpha_0 - \beta| < \varepsilon$. But then $\hat{f}(z)$ extends to a holomorphic function in an open neighborhood of $\alpha_0$, and, therefore, $\limsup_{n \to \infty} C(f,\alpha_0,n) < \infty$. This implies that the sequence $C(f,\alpha_0,n)$ is bounded and $f \notin \bigcap_{N \in \mathbb{N}} \Omega(\alpha_0,N)$, since $C(f,\alpha_0,n) \leq N_0$, for every $n \in \mathbb{N}$, when the integer $N_0$ is chosen so that $N_0 > \sup_{n \in \mathbb{N}} C(f,\alpha_0,n)$. This contradicts our assumption that $f \in \bigcap_{\alpha \in A} \bigcap_{N \in \mathbb{N}} \Omega(\alpha,N)$ and completes the proof of the claim.

Finally, by Baire's theorem and Lemma 3, the set $\bigcap_{\alpha \in A} \bigcap_{N \in \mathbb{N}} \Omega(\alpha,N)$ is dense and $\mathcal{G}_\delta$ in $\mathcal{S}([0,\infty))$, and this completes the proof of the theorem.

## 3. Generalizations

We can also state and prove a more general version of Theorem 1, by setting

$$\rho_{\ell,k}(f) := \sup_{t \geq 0}\left[ t^\ell \left|\frac{\partial^k f(t)}{\partial t^k}\right|\right]$$

and considering the space

$$\mathcal{S}^m([0,\infty)) := \{f \in C^m([0,\infty)) : \rho_{\ell,k}(f) < \infty \text{ for every nonnegative integers } \ell \text{ and } k \text{ with } k \leq m\},$$

for each fixed $m = 0, 1, 2, \ldots, \infty$. ($\mathcal{S}^\infty([0,\infty))$ is the space $\mathcal{S}([0,\infty))$ which we considered previously.)
The space $\mathcal{S}^m([0,\infty))$ becomes a Frechet space with the topology induced by the seminorms $\rho_{\ell,k}$, $k \leq m$.

The proof of the following theorem is similar to the proof of Theorem 1.

**Theorem 2** *The set $\{f \in \mathcal{S}^m([0,\infty)) : r(f,\alpha) = \infty \text{ for every } \alpha \in \mathbb{R}\}$ is dense and $\mathcal{G}_\delta$ in the space $\mathcal{S}^m([0,\infty))$.*

We also give a generalization of the example which we considered Section 2.

For $t \geq 0$, we consider the function $\psi_p(t) := e^{-t^p}$. ($\chi_0$ is the case $p = 1/2$ of $\psi_p$.)

**Proposition** *Let $p$ be a real number with $0 < p < 1$. Then the function*

$$\hat{\psi}_p(z) = \int_0^\infty e^{-izt} \psi_p(t) dt,$$

*which is defined – a priori – for $z \in \mathbb{C}$ with $\operatorname{Im} z \leq 0$, extends analytically to $z \in \mathbb{C} - i(\mathbb{R}^+ \cup \{0\})$ and $0$ is singular point of this function in the sense that it does not extend analytically to any open neighborhood of $0$.*

*Proof* For $z \in \mathbb{C}$ with $\operatorname{Im} z \leq 0$, we have

$$\hat{\psi}_p^{(n)}(z) = \int_0^\infty (-it)^n e^{-izt} e^{-t^p} dt.$$

(As before, when $z \in \mathbb{R}$, $\hat{\psi}_p^{(n)}(z)$, $n = 0, 1, 2, \ldots$, are the real derivatives of the restriction of the function $\hat{\psi}_p$ to the real line.)
It follows (see [74, p. 364]) that

$$\hat{\psi}_p^{(n)}(0) = \int_0^\infty (-it)^n e^{-t^p} dt = (-i)^n \frac{1}{p} \Gamma\left(\frac{n+1}{p}\right) \approx (-i)^n \sqrt{2\pi\left(\frac{n+1}{p} - 1\right)} \left[\frac{1}{e}\left(\frac{n+1}{p} - 1\right)\right]^{\left(\frac{n+1}{p} - 1\right)}, \text{ as } n \to \infty.$$

Also $n! \approx \sqrt{2\pi n}(n/e)^n$, and, therefore, the radius of convergence of the power series is $\sum_{n=0}^\infty [\hat{\psi}_p^{(n)}(0)/n!] z^n$ is zero.
This proves that, indeed, $0$ is singular point of $\hat{\psi}_p(z)$.

In order to obtain the analytic continuation of $\hat{\psi}_p(z)$ to $\mathbb{C} - i(\mathbb{R}^+ \cup \{0\})$, we comptute $\hat{\psi}_p(-i\mu)$, for $\mu > 0$. We have

$$\hat{\psi}_p(-i\mu) = \int_0^\infty e^{-\mu t} e^{-t^p} dt = p\int_0^\infty e^{-\mu s^{1/p}} e^{-s} s^{p-1} dt$$

$$= p\int_0^\infty e^{-\mu s^{1/p}} \sum_{n=0}^\infty \frac{(-1)^n}{n!} s^n s^{p-1} dt = p\sum_{n=0}^\infty \frac{(-1)^n}{n!} \int_0^\infty e^{-\mu s^{1/p}} s^{n+p-1} ds$$

$$= \frac{p^2}{\mu^{p^2}} \sum_{n=0}^{\infty} \frac{(-1)^n}{n!} \frac{1}{\mu^{np}} \Gamma((n+p)p).$$

Indeed, since

$$\Gamma((n+p)p) \approx \sqrt{2\pi(np+p^2-1)} \left[\frac{1}{e}(np+p^2-1)\right]^{np+p^2-1}, \text{ as } n \to \infty,$$

we obtain that

$$\int_0^{\infty} \sum_{n=0}^{\infty} \left| e^{-\mu s^{1/p}} \frac{(-1)^n}{n!} s^n s^{p-1} \right| dt = \frac{p}{\mu^{p^2}} \sum_{n=0}^{\infty} \frac{1}{n!} \frac{1}{\mu^{np}} \Gamma((n+p)p) < \infty,$$

and this justifies the interchange of the integration and the summation, that we used previously. Moreover, setting $z = -i\mu$, i.e., $\mu = iz$, we obtain the required analytic continuation of $\hat{\psi}_p(z)$ to $\mathbb{C} - i(\mathbb{R}^+ \cup \{0\})$, given by the series

$$\frac{p^2}{(iz)^{p^2}} \sum_{n=0}^{\infty} \frac{(-1)^n}{n!} \frac{1}{(iz)^{np}} \Gamma((n+p)p).$$

It remains to point out that we choose the branches of the power functions $(iz)^{p^2}$ and $(iz)^{np}$ to be the ones defined from the logarithm $\log w = \log|w| + i \arg w$, with $-\pi < \arg w \le \pi$, so that the functions

$$(iz)^{p^2} = \exp[p^2 \log(iz)] \text{ and } (iz)^{np} = \exp[np \log(iz)]$$

are holomorfphic for $z$ in $\mathbb{C} - i(\mathbb{R}^+ \cup \{0\})$.

To state the next theorem, we define

$$\mathrm{K}_M(f,\alpha,n) = \frac{1}{n^M} \sqrt[n]{\frac{|(\hat{f})^{(n)}(\alpha)|}{n!}},$$

for $f \in \mathcal{S}^0([0,\infty))$, $\alpha \in \mathbb{R}$, $n \in \mathbb{N}$ and $M \in \mathbb{N} \cup \{0\}$. (The case $\mathrm{K}_0(f,\alpha,n) = C(f,\alpha,n)$.)

**Theorem 3** *Let* $\mathrm{A} \subset \mathbb{R}$ *be a countable dense subset in the real line. Then the set*

$$\mathcal{B} := \{f \in \mathcal{S}^0([0,\infty)) : \limsup_{n \to \infty} \mathrm{K}_M(f,\alpha,n) = \infty \text{ for every } M \in \mathbb{N} \text{ and } \alpha \in \mathrm{A}\}$$

*is dense and* $\mathcal{G}_\delta$ *in the space* $\mathcal{S}^0([0,\infty))$.

**Proof** For fixed $N \in \mathbb{N}$, $M \in \mathbb{N}$ and $\alpha \in \mathrm{A}$, we consider the set

$$\Theta_M(\alpha,N) := \bigcup_{n \in \mathbb{N}} \{f \in \mathcal{S}^0([0,\infty)) : \mathrm{K}_M(f,\alpha,n) > N\}.$$

Then

$$\bigcap_{M \in \mathbb{N}} \bigcap_{\alpha \in \mathrm{A}} \bigcap_{N \in \mathbb{N}} \Theta_M(\alpha,N) \subseteq \mathcal{B}.$$

Indeed, if $f \in \mathcal{S}^0([0,\infty))$ and $f \notin \mathcal{B}$ then, for some $M_0 \in \mathbb{N}$ and $\alpha_0 \in \mathrm{A}$, $\limsup_{n \to \infty} \mathrm{K}_{M_0}(f,\alpha_0,n) < \infty$, i.e., the sequence $\mathrm{K}_{M_0}(f,\alpha_0,n)$ is bounded. But then we can choose $N_0$ so that $\mathrm{K}_{M_0}(f,\alpha_0,n) < N_0$, for every $n \in \mathbb{N}$, and this imlies that $f \notin \Theta_{M_0}(\alpha_0,N_0)$.

Since, clearly, we also have $\bigcap_{M \in \mathbb{N}} \bigcap_{\alpha \in \mathrm{A}} \bigcap_{N \in \mathbb{N}} \Theta_M(\alpha,N) \subseteq \mathcal{B}$, it follows that $\mathcal{B} = \bigcap_{M \in \mathbb{N}} \bigcap_{\alpha \in \mathrm{A}} \bigcap_{N \in \mathbb{N}} \Theta_M(\alpha,N)$.

Thus, in view of Baire's theorem, it suffices to show that, for fixed $N \in \mathbb{N}$, $M \in \mathbb{N}$ and $\alpha \in A$, the set $\Theta_M(\alpha, N)$ open and dense in the space $\mathcal{S}^0([0,\infty))$.

First, for fixed $n$, $N$, $M$ and $\alpha$,

$$K_M(f,\alpha,n) \leq N \Leftrightarrow \left|(\hat{f})^{(n)}(\alpha)\right| \leq N^n n^{Mn} n! \Leftrightarrow \left|\int_0^\infty e^{-i\alpha t}(-it)^n f(t)dt\right| \leq N^n n^{Mn} n!. \tag{3}$$

This implies that the set $\{f \in \mathcal{S}^0([0,\infty)) : K_M(f,\alpha,n) \leq N\}$ is closed in $\mathcal{S}^0([0,\infty))$, and, therefore, $\Theta_M(\alpha, N)$ is indeed open.

It remains to show that $\Theta_M(\alpha, N)$ is dense in $\mathcal{S}^0([0,\infty))$. For this purpose, let us fix a $f \in \mathcal{S}^0([0,\infty)) - \Theta_M(\alpha, N)$. Then, in view of (3),

$$\left|(\hat{f})^{(n)}(\alpha)\right| \leq N^n n^{Mn} n!, \text{ for every } n.$$

Now, we consider the function

$$X_{\alpha,M}(t) := e^{i\alpha t} \exp\left(-t^{1/M^2}\right), \ t \geq 0.$$

Then

$$[(X_{\alpha,M})\hat{\ }]^{(n)}(z) = \int_0^\infty e^{-izt}(-it)^n e^{i\alpha t} \exp(-t^{1/M^2})dt$$

and, therefore,

$$[(X_{\alpha,M})\hat{\ }]^{(n)}(\alpha) = \int_0^\infty (-it)^n \exp\left(-t^{1/M^2}\right)dt$$

$$= (-i)^n M^2 \int_0^\infty s^{M^2 n} s^{M^2-1} \exp(-s)ds = (-i)^n M^2 [(M^2 n + M^2 - 1)!].$$

We claim that

$$f + \frac{1}{j} X_{\alpha,M} \in \Theta_M(\alpha, N), \text{ for every } j \in \mathbb{N}. \tag{4}$$

In order to prove this, it suffices to show that, for fixed $j$, there is $n_0 \in \mathbb{N}$ so that

$$K_M(f + \tfrac{1}{j} X_{\alpha,M}, \alpha, n) > N, \text{ i.e., } \left|\hat{f}^{(n)}(\alpha) + \frac{1}{j}[(X_{\alpha,M})\hat{\ }]^{(n)}(\alpha)\right| > N^n n^{Mn} n!, \text{ for } n \geq n_0.$$

Since

$$\left|\hat{f}^{(n)}(\alpha) + \frac{1}{j}[(X_{\alpha,M})\hat{\ }]^{(n)}(\alpha)\right| \geq \frac{1}{j}\left|[(X_{\alpha,M})\hat{\ }]^{(n)}(\alpha)\right| - \left|\hat{f}^{(n)}(\alpha)\right| \text{ and } \left|(\hat{f})^{(n)}(\alpha)\right| \leq N^n n^{Mn} n!,$$

it suffices to show that, given $M$, $N$ and $m$, there is $n_0 \in \mathbb{N}$ so that

$$M^2[(M^2 n + M^2 - 1)!] \geq 2jN^n n^{Mn} n!, \text{ for } n \geq n_0.$$

But, with $n \geq N$ and $M \geq 2$, we have $M^2 n + M^2 - 1 - 2n \geq Mn + 1$, and, therefore,

$$(M^2 n + M^2 - 1)! \geq n! n^{M^2 n + M^2 - 1 - n} \geq n! N^n n^{M^2 n + M^2 - 1 - 2n} \geq nN^n n^{Mn} n!.$$

This proves (4) and completes the proof of the theorem.

**Remark** Since $n! \approx \sqrt{2\pi n}(n/e)^n$ implies $\sqrt[n]{n!} \approx n/e$, we have

$$\limsup_{n \to \infty} K_M(f,\alpha,n) = e \limsup_{n \to \infty} \frac{1}{n^{M+1}} \sqrt[n]{\left|(\hat{f})^{(n)}(\alpha)\right|}.$$

Therefore, for a fixed $f \in \mathcal{S}^0([0,\infty))$,

$$\limsup_{n\to\infty} K_M(f,\alpha,n) = \infty \text{ for every } M \in \mathbb{N} \Leftrightarrow \limsup_{n\to\infty} \frac{1}{n^M}\sqrt[n]{|(\hat{f})^{(n)}(\alpha)|} = \infty \text{ for every } M \in \mathbb{N}.$$

**Comments** 1 We leave open the question whether, for a fixed $m \geq 1$, the set

$$\{f \in \mathcal{S}^m([0,\infty)) : \limsup_{n\to\infty} K_M(f,\alpha,n) = \infty \text{ for every } M \in \mathbb{N} \text{ and for every } \alpha \in A\}$$

is dense and $\mathcal{G}_\delta$ in the space $\mathcal{S}^m([0,\infty))$.

**2** We also leave open the question whether the set

$$\{f \in \mathcal{S}^0([0,\infty)) : \limsup_{n\to\infty} K_M(f,\alpha,n) = \infty \text{ for every } M \in \mathbb{N} \text{ and for every } \alpha \in \mathbb{R}\}$$

is dense and $\mathcal{G}_\delta$ in the space $\mathcal{S}^0([0,\infty))$.